# Finite q-identities related to well-known theorems of Euler and Gauss

Johann Cigler

Fakultät für Mathematik
Universität Wien
A-1090 Wien, Nordbergstraße 15

email: johann.cigler@univie.ac.at

**Abstract**
*We give generalizations of a finite version of Euler's pentagonal number theorem and of a $q-$identity of Gauss by introducing a new parameter.*

## 0. Introduction

Let $\begin{bmatrix} n \\ k \end{bmatrix} = \frac{(1-q^{n-k+1})\cdots(1-q^n)}{(1-q)\cdots(1-q^k)}$ be a $q-$binomial coefficient.

A. Berkovich and F. G. Garvan [1] and with another method S.O. Warnaar [4] have proved that

$$\sum_{j=-L}^{2L} (-1)^j q^{\frac{j(3j+1)}{2}} \begin{bmatrix} 2L-j \\ L+j \end{bmatrix} = 1$$

for all $L \in \mathbb{N}$. This is a finite version of Euler's pentagonal number theorem, because for $L \to \infty$ it reduces to the pentagonal number theorem

$$\prod_{n=1}^{\infty}(1-q^n) = \sum_{j=-\infty}^{\infty} (-1)^j q^{\frac{j(3j+1)}{2}}.$$

Let $r_n(x,a) = \sum_k \begin{bmatrix} n \\ k \end{bmatrix} x^k a^{n-k}$. A famous theorem of Gauss states that $r_{2n+1}(1,-1) = 0$ and

$$r_{2n}(1,-1) = \frac{(q;q)_{2n}}{(q^2;q^2)_n} = (1-q)(1-q^2)\cdots(1-q^{2n-1}).$$

Generalizing these results we give explicit evaluations of the sums

$$h(L,k) = \sum_{j=-L}^{2L} (-1)^j q^{\frac{j(3j+1)}{2}+kj} \begin{bmatrix} 2L-j \\ L+j \end{bmatrix}$$

and

$$r_n(1,-q^k) = \sum_{j=0}^{n} \begin{bmatrix} n \\ j \end{bmatrix} (-1)^j q^{kj}.$$

# 1. Variations of a finite version of Euler's pentagonal number theorem

Generalizing the formula of Berkovich and Garvan we give explicit evaluations of the sums

$$h(L,k) = \sum_{j=-L}^{2L} (-1)^j q^{\frac{j(3j+1)}{2}+kj} \begin{bmatrix} 2L-j \\ L+j \end{bmatrix}. \tag{1.1}$$

For small values of $k \in \mathbb{N}$ we get the following formulas for $h(L,k)$:

$h(L,0) = 1$

$h(L,1) = q^L \begin{bmatrix} 1 \\ 1 \end{bmatrix}$

$h(L,2) = -q^{-1} \begin{bmatrix} 2 \\ 0 \end{bmatrix} + q^{L-1} \begin{bmatrix} 2 \\ 1 \end{bmatrix}$

$h(L,3) = -q^{-2} \begin{bmatrix} 3 \\ 0 \end{bmatrix} + q^{2L-2} \begin{bmatrix} 3 \\ 2 \end{bmatrix} - q^{3L-1} \begin{bmatrix} 3 \\ 3 \end{bmatrix}$

$h(L,4) = -q^{L-4} \begin{bmatrix} 4 \\ 1 \end{bmatrix} + q^{2L-4} \begin{bmatrix} 4 \\ 2 \end{bmatrix} - q^{4L-2} \begin{bmatrix} 4 \\ 4 \end{bmatrix}$

For $q=1$ the left hand side of each equation reduces to 1. The right hand side reduces to a difference of sums of the form $A_{k,i} = \sum_{j \equiv i \,(\text{mod}\, 3)} \binom{k}{j}$. E.g.

$A_{7,1} - A_{7,2} = \left( \binom{7}{1} + \binom{7}{4} + \binom{7}{7} \right) - \left( \binom{7}{2} + \binom{7}{5} \right) = (7+35+1) - (21+21) = 1$. The existence of such a representation becomes obvious from

$(1+\rho)^n = A_{n,0} + A_{n,1}\rho + A_{n,2}\rho^2 = (-\rho^2)^n$,

where $\rho$ denotes a third root of unity.

E.g. if $n = 3m+1$ we have $(-\rho^2)^{3m+1} = (-1)^{m+1} \rho^2 = (-1)^m (1+\rho)$. Therefore $A_{n,1} - A_{n,2} = (-1)^m$.

In the general case we define a sequence $w(n)$ for $n \in \mathbb{Z}$ by $w(n) = 0$ if $n \equiv 2 (\text{mod}\, 3)$ and $w(n) = (-1)^{\lfloor \frac{n}{3} \rfloor} q^{-\frac{n(n-1)}{6}}$ else. This means

$$w(3k) = (-1)^k q^{-\frac{k(3k-1)}{2}}$$
$$w(3k+1) = (-1)^k q^{-\frac{k(3k+1)}{2}} \tag{1.2}$$
$$w(3k+2) = 0.$$

From the definition follows that
$$w(-n) = w(n+1). \tag{1.3}$$

Then we have

**Theorem 1**

For each $L \in \mathbb{N}$ and $k \in \mathbb{N}$ the sum $h(L,k) = \sum_{j=-L}^{2L} (-1)^j q^{\frac{j(3j+1)}{2}+kj} \begin{bmatrix} 2L-j \\ L+j \end{bmatrix}$ has the value

$$h(L,k) = \sum_{j=0}^{k} q^{\binom{j+1}{2}} q^{jL} \begin{bmatrix} k \\ j \end{bmatrix} (-1)^j w(-k-j) = \sum_{j=0}^{k} q^{\binom{j+1}{2}} q^{jL} \begin{bmatrix} k \\ j \end{bmatrix} (-1)^j w(k+j+1). \tag{1.4}$$

Since $w(-k-j) = 0$ if $-k-j \equiv 2 \pmod 3$ there remain only $q$-binomial coefficients $\begin{bmatrix} k \\ j \end{bmatrix}$ with $j \equiv -k \pmod 3$ or $j \equiv -k+2 \pmod 3$. The sign of the coefficient of $q^{jL}$ is given by $(-1)^{j+\left\lfloor -\frac{k+j}{3} \right\rfloor} = (-1)^{j+3+\left\lfloor -\frac{k+j+3}{3} \right\rfloor}$. Therefore all terms with $j$ in the same residue class mod 3 have the same sign.

E.g. $h(L,7) = A_{7,1}(L,q) - A_{7,2}(L,q)$ with $A_{7,1}(L,q) = q^{L-11} \begin{bmatrix} 7 \\ 1 \end{bmatrix} + q^{4L-12} \begin{bmatrix} 7 \\ 4 \end{bmatrix} + q^{7L-7} \begin{bmatrix} 7 \\ 7 \end{bmatrix}$ and $A_{7,2}(L,q) = q^{2L-12} \begin{bmatrix} 7 \\ 2 \end{bmatrix} + q^{5L-11} \begin{bmatrix} 7 \\ 5 \end{bmatrix}$.

The expansion of $h(L,n)$ has a vanishing constant term if and only if $n \equiv 1 \pmod 3$. Therefore we get $h(3k+1) = \lim_{L \to \infty} h(L, 3k+1) = 0$.

For $n = 3k$ the constant term is $w(-3k) = (-1)^k q^{-\frac{k(3k+1)}{2}}$ and therefore $h(3k) = (-1)^k q^{-\frac{k(3k+1)}{2}}$.
Analogously we get $h(3k-1) = (-1)^k q^{-\frac{k(3k-1)}{2}}$.
From the definition of $h(L,n)$ we see that
$$h(3k+i) \prod_{n=1}^{\infty} (1-q^n) = \sum_{j=-\infty}^{\infty} (-1)^j q^{\frac{j(3j+1)}{2}+(3k+i)j}.$$
This can also be verified by using Jacobi's triple product identity. The case $h(0) = 1$ gives Euler's pentagonal number theorem.
For the triple product implies

$$\sum_{j \in \mathbb{Z}} (-1)^j q^{aj^2+bj} = \prod_{n=0}^{\infty} (1-q^{2an+a+b})(1-q^{2an+a-b})(1-q^{2an+2a}).$$

Therefore
$$\sum_{j=-\infty}^{\infty} (-1)^j q^{\frac{j(3j+1)}{2}+(3k+i)j} = \prod_{n \geq 0} (1-q^{3n+2+3k+i})(1-q^{3n+1-3k-i})(1-q^{3n+3}).$$

For $i = 3l+1$ the product vanishes because one of the middle factors vanishes.
For $i = 3l$ we get

$$\prod_{n\geq 0}(1-q^{3n+2+3k+i})(1-q^{3n+1-3k-i})(1-q^{3n+3}) = q^{-\frac{(k+l)(3(k+l)+1)}{2}}\prod_{n\geq 1}(1-q^n) = h(3k+i)\prod_{n\geq 1}(1-q^n)$$

because to each term $1-q^{3n+2}, 0 \leq n < k+l$, which does not occur in the first factors of the product there corresponds precisely one term $1-q^{-(3n+2)}$ in the middle factors and $2+5+\cdots+(3(k+l)-1) = \frac{(k+l)(3(k+l)+1)}{2}$. The case $i = 3l-1$ can be treated in the same way.

**Proof of Theorem 1**

In order to prove this theorem we use some results about the $q$-Fibonacci polynomials

$$F_n(x,s) = \sum_{k=0}^{\lfloor\frac{n-1}{2}\rfloor}\begin{bmatrix}n-k-1\\k\end{bmatrix}q^{\binom{k+1}{2}}x^{n-1-2k}s^k$$ from our previous paper [2]. These are $q$-analogues

of the Fibonacci polynomials $f_n(x,s) = \sum_{k=0}^{\lfloor\frac{n-1}{2}\rfloor}\binom{n-1-k}{k}x^{n+2k-1}s^k$, which are characterized by the recurrence relation $f_n(x,s) = xf_{n-1}(x,s) + sf_{n-2}(x,s)$ and the initial values $f_0(x,s) = 0$ and $f_1(x,s) = 1$.

Let

$$G(L,i,s) = \sum_{j=-L}^{2L+i} s^j q^{\frac{j(3j-1)}{2}-ij}\begin{bmatrix}2L+i-j\\L+j\end{bmatrix} \quad (1.5)$$

and

$$f_n(s) = \sum_{k=0}^{n-1}\begin{bmatrix}n-1-k\\k\end{bmatrix}_{\frac{1}{q}} q^{-\binom{k+1}{2}}s^k. \quad (1.6)$$

Then $f_n(s) = F_n(1,s)\big|_{q\to\frac{1}{q}}$ and therefore

$$f_n(qs) = f_{n-1}(s) + sf_{n-2}(s) \quad (1.7)$$

by [2] (2.2).

Using the easily verified formula $q^{-\binom{k}{2}} \begin{bmatrix} n-k \\ k \end{bmatrix}_{\frac{1}{q}} = q^{-\binom{n}{2}} q^{k^2+\binom{n-k}{2}} \begin{bmatrix} n-k \\ k \end{bmatrix}$

we get

$$\sum_{j=-L}^{2L+i} s^{j+L} q^{\frac{j(3j-1)}{2} - ij} \begin{bmatrix} 2L+i-j \\ L+j \end{bmatrix} = q^{\frac{L(3L+1)}{2} + iL} \sum_{k=0}^{3L+i} q^{-\binom{k+1}{2}} \begin{bmatrix} 3L+i-k \\ k \end{bmatrix}_{\frac{1}{q}} s^k. \tag{1.8}$$

Then (1.8) can be written as

$$G(L,i,s) = q^{\frac{L(3L+1)}{2}+iL} s^{-L} f_{3L+i+1}(s). \tag{1.9}$$

Combining (1.7) and (1.9) we get
$$G(L,i,qs) = G(L,i-1,s) + q^L s G(L,i-2,s). \tag{1.10}$$

By [2] (3.2)
$$f_{3n}(-q) = 0, \quad f_{3n+1}(-q) = (-1)^n q^{-\frac{n(3n-1)}{2}}, \quad f_{3n+2}(-q) = (-1)^n q^{-\frac{n(3n+1)}{2}}. \tag{1.11}$$

The formula $F_{-n}(1,s) = (-1)^{n-1} \frac{F_n(1,s)}{s^n}$ (cf. [2] (2.7))

gives

$F_{-n}(1,-\frac{1}{q}) = -F_n(1,-\frac{1}{q}) q^n.$

Therefore we get

$F_{-3k+1}(1,-\frac{1}{q}) = (-1)^k q^{\frac{k(3k+1)}{2}}, \quad F_{-3k+2}(1,-\frac{1}{q}) = (-1)^k q^{\frac{k(3k-1)}{2}}, \quad F_{-3k}(1,-\frac{1}{q}) = 0.$

This implies that (1.11) holds for all $n \in \mathbb{Z}$.

We use (1.9) to extend $G(L,i,s)$ to values with $3L+i+1 < 0$.
Then it is easy to verify that

$$G(L,n,-q) = w(n) \tag{1.12}$$

for all $n \in \mathbb{Z}$. Note that the right-hand side does not depend on $L \in \mathbb{N}$.

Formula (1.10) gives immediately

$$G(L,i,q^k s) = \sum_{j=0}^{k} q^{\binom{j}{2}} q^{jL} s^j \begin{bmatrix} k \\ j \end{bmatrix} G(L,i-k-j,s). \qquad (1.13)$$

This follows by induction from

$$G(L,i,q^{k+1}s) = G(L,i-1,q^k s) + q^{L+k} s G(L,i-2,q^k s)$$

$$= \sum_{j=0}^{k} q^{\binom{j}{2}} q^{jL} s^j \begin{bmatrix} k \\ j \end{bmatrix} G(L,i-k-j-1,s) + \sum_{j=0}^{k} q^{\binom{j}{2}} q^{jL+L} q^k s^j \begin{bmatrix} k \\ j \end{bmatrix} G(L,i-k-j-2,s)$$

$$= \sum_{j \geq 0} q^{\binom{j}{2}} q^{jL} s^j \begin{bmatrix} k \\ j \end{bmatrix} G(L,i-k-j-1,s) + \sum_{j \geq 0} q^{\binom{j}{2}} q^{jL} s^j q^k \begin{bmatrix} k \\ j \end{bmatrix} G(L,i-k-j-1,s)$$

$$= \sum_{j=0}^{k+1} q^{\binom{j}{2}} q^{jL} s^j \begin{bmatrix} k+1 \\ j \end{bmatrix} G(L,i-k-1-j,s).$$

The theorem follows if in (1.13) we set $s = -q$ and $i = 0$.

## 2. Variations of a q-identity of Gauss

Consider now the Rogers-Szegö polynomials

$$r_n(x,a) = \sum_{k} \begin{bmatrix} n \\ k \end{bmatrix} x^k a^{n-k}. \qquad (2.1)$$

Gauss's theorem states that
$r_{2n+1}(1,-1) = 0$ and

$$r_{2n}(1,-1) = \frac{(q;q)_{2n}}{(q^2;q^2)_n} = (1-q)(1-q^2)\cdots(1-q^{2n-1}). \qquad (2.2)$$

A very simple proof uses Euler's $q-$exponential series

$$e(x) = \frac{1}{(x;q)_\infty} = \frac{1}{(1-x)(1-qx)(1-q^2 x)\cdots} = \sum_{k \geq 0} \frac{x^k}{(q;q)_k}. \qquad (2.3)$$

The generating function of the Rogers-Szegö polynomials is given by

$$\sum_n \frac{r_n(x,a)}{(q;q)_n} z^n = e(xz)e(az). \qquad (2.4)$$

This follows from

$$e(xz)e(az) = \sum_k \frac{x^k z^k}{(q;q)_k} \sum_\ell \frac{a^\ell z^\ell}{(q;q)_\ell} = \sum_{k,\ell} \frac{z^{k+\ell}}{(q;q)_{k+\ell}} \frac{(q;q)_{k+\ell}}{(q;q)_k (q;q)_\ell} x^k a^\ell = \sum_n \frac{z^n}{(q;q)_n} \sum_{k+\ell=n} \begin{bmatrix} k+\ell \\ k \end{bmatrix} x^k a^\ell$$

$$= \sum_n \frac{r_n(x,a)}{(q;q)_n} z^n.$$

As a special case we get

$$\sum_n \frac{r_n(1,-1)}{(q;q)_n} z^n = e(-z)e(z) = \frac{1}{(1+z)(1+qz)(1+q^2z)\cdots} \frac{1}{(1-z)(1-qz)(1-q^2z)\cdots}$$

$$= \frac{1}{(1-z^2)(1-q^2z^2)(1-q^4z^2)\cdots} = \sum_n \frac{(z^2)^n}{(q^2;q^2)_n} = \sum_n \frac{z^{2n}}{(q;q)_{2n}} (1-q)(1-q^3)\cdots(1-q^{2n-1}),$$
(2.5)

which is equivalent with Gauss's theorem by comparing coefficients.

B.A Kupershmidt [3] has given a formula for $r_n(1,-x) = (-1)^n r_n(x,-1)$ which can be easily deduced from the generating function

$$\sum_n \frac{r_n(x,-1)}{(q;q)_n} z^n = e(xz)e(-z) = \left(\frac{e(xz)}{e(z)}\right)(e(z)e(-z)).$$

He raised the problem of explicitly evaluating $r_n(1,-q^k)$. In the following we give two such formulas.

For the first one we generalize the method we used to prove Gauss's theorem:

From $e(qx) = \dfrac{1}{(1-qx)(1-q^2x)\cdots} = (1-x)e(x)$

we get

$$\sum_n \frac{r_n(1,-q)}{(q;q)_n} z^n = e(-qz)e(z) = \frac{1+z}{(1+z)(1+qz)(1+q^2z)\cdots} \frac{1}{(1-z)(1-qz)(1-q^2z)\cdots}$$

$$= \frac{1+z}{(1-z^2)(1-q^2z^2)(1-q^4z^2)\cdots} = (1+z)\sum_n \frac{(z^2)^n}{(q^2;q^2)_n}.$$

This implies $r_{2n}(1,-q) = r_{2n}(1,-1)$ and $r_{2n-1}(1,-q) = r_{2n}(1,-1)$.

More generally for $k \geq 1$

$$\sum_n \frac{r_n(1,-q^k)}{(q;q)_n} z^n = e(-q^k z)e(z) = \frac{(1+z)\cdots(1+q^{k-1}z)}{(1+z)(1+qz)(1+q^2z)\cdots} \frac{1}{(1-z)(1-qz)(1-q^2z)\cdots}$$

$$= (1+z)\cdots(1+q^{k-1}z)\sum_n \frac{r_n(1,-1)}{(q;q)_n} z^n = \sum_{j=0}^k q^{\binom{j}{2}} \begin{bmatrix} k \\ j \end{bmatrix} z^j \sum_n \frac{z^{2n}}{(q^2;q^2)_n}.$$

Comparing coefficients we get

$$r_n(1,-q^k) = (q;q)_n \sum_{j+2\ell=n} q^{\binom{j}{2}} \begin{bmatrix} k \\ j \end{bmatrix} \frac{1}{(q^2;q^2)_\ell}.$$
(2.6)

This may be written in the form

$$r_{2n-1}(1,-q^k) = r_{2n}(1,-1) \sum_{0 \leq \ell < n} q^{\binom{2n-2\ell-1}{2}} \begin{bmatrix} k \\ 2n-2\ell-1 \end{bmatrix} \frac{(q^2;q^2)_{n-1}}{(q^2;q^2)_\ell}$$

$$= r_{2n}(1,-1) \sum_{j \geq 0} q^{\binom{2j+1}{2}} \begin{bmatrix} k \\ 2j+1 \end{bmatrix} \prod_{i=n-j}^{n-1} (1-q^{2i})$$

and

$$r_{2n}(1,-q^k) = r_{2n}(1,-1)\sum_{j\geq 0} q^{\binom{2j}{2}} \begin{bmatrix} k \\ 2j \end{bmatrix} \prod_{i=n-j+1}^{n}(1-q^{2i}).$$

The first values are

$$\frac{r_{2n-1}(1,-q)}{r_{2n}(1,-1)} = 1$$

$$\frac{r_{2n-1}(1,-q^2)}{r_{2n}(1,-1)} = \begin{bmatrix} 2 \\ 1 \end{bmatrix}$$

$$\frac{r_{2n-1}(1,-q^3)}{r_{2n}(1,-1)} = \begin{bmatrix} 3 \\ 1 \end{bmatrix} + q^3(1-q^{2n-2})$$

$$\frac{r_{2n-1}(1,-q^4)}{r_{2n}(1,-1)} = \begin{bmatrix} 4 \\ 1 \end{bmatrix} + q^3 \begin{bmatrix} 4 \\ 3 \end{bmatrix}(1-q^{2n-2})$$

$$\frac{r_{2n-1}(1,-q^5)}{r_{2n}(1,-1)} = \begin{bmatrix} 5 \\ 1 \end{bmatrix} + \begin{bmatrix} 5 \\ 3 \end{bmatrix} q^3(1-q^{2n-2}) + \begin{bmatrix} 5 \\ 5 \end{bmatrix} q^{10}(1-q^{2n-2})(1-q^{2n-4})$$

and

$$\frac{r_{2n}(1,-q)}{r_{2n}(1,-1)} = 1$$

$$\frac{r_{2n}(1,-q^2)}{r_{2n}(1,-1)} = 1 + q(1-q^{2n})$$

$$\frac{r_{2n}(1,-q^3)}{r_{2n}(1,-1)} = 1 + q(1-q^{2n})\begin{bmatrix} 3 \\ 2 \end{bmatrix}$$

$$\frac{r_{2n}(1,-q^4)}{r_{2n}(1,-1)} = 1 + q(1-q^{2n})\begin{bmatrix} 4 \\ 2 \end{bmatrix} + q^6(1-q^{2n})(1-q^{2n-2})\begin{bmatrix} 4 \\ 4 \end{bmatrix}.$$

These are polynomials in $q^{2n}$. Now we want to compute the coefficients of these polynomials.
In order to do this we start from the formula

$$r_n(q^2 x, a) - (1 - \frac{qx}{a})r_n(qx, a) - q^{n+1}\frac{x}{a}r_n(x,a) = 0, \qquad (2.7)$$

which is easily verified by comparing coefficients:

$$q^{2k}\begin{bmatrix} n \\ k \end{bmatrix} - q^k \begin{bmatrix} n \\ k \end{bmatrix} + q^k \begin{bmatrix} n \\ k-1 \end{bmatrix} - q^{n+1}\begin{bmatrix} n \\ k-1 \end{bmatrix} = q^k\left((q^k-1)\begin{bmatrix} n \\ k \end{bmatrix} + (1-q^{n+1-k})\begin{bmatrix} n \\ k-1 \end{bmatrix}\right) = 0.$$

Let now
$$b(n,k) = -\frac{r_{2n-1}(q^k,-1)}{r_{2n}(1,-1)}. \tag{2.8}$$

Then we get
$$b(n,k+2) - (1+q^{k+1})b(n,k+1) + q^{2n+k}b(n,k) = 0. \tag{2.9}$$

Define now polynomials $f(k,s)$ by the recurrence
$$f(k,s) = (1+q^{k-1})f(k-1,s) - q^{k-2}sf(k-2,s) \tag{2.10}$$

and initial values $f(0,s) = 0$ and $f(1,s) = 1$.
Then
$$b(n,k) = f(k, q^{2n}). \tag{2.11}$$

The polynomial $f(k,s)$ is a $q$-analogue of the Fibonacci polynomial
$$f_k(2,-s) = \sum_{j=0}^{\lfloor \frac{k-1}{2} \rfloor} (-1)^j s^j \binom{k-1-j}{j} 2^{k-1-2j}.$$

It is easily verified that this formula has the direct $q$-analogue
$$f(k,s) = \sum_{j=0}^{\lfloor \frac{k-1}{2} \rfloor} (-1)^j q^{j^2} s^j \begin{bmatrix} k-j-1 \\ j \end{bmatrix}_{q^2} \prod_{i=1}^{k-1-2j}(1+q^i). \tag{2.12}$$

For by comparing coefficients the recursion (2.10) is equivalent with the identity
$$\begin{bmatrix} k-j-1 \\ j \end{bmatrix}_{q^2}(1+q^{k-1-2j}) - (1+q^{k-1})\begin{bmatrix} k-j-2 \\ j \end{bmatrix}_{q^2} - q^{k-1-2j}\begin{bmatrix} k-j-2 \\ j-1 \end{bmatrix}_{q^2}(1+q^{k-1-2j}) = 0.$$

This is trivial, because we get
$$\begin{bmatrix} k-j-1 \\ j \end{bmatrix}_{q^2} - \begin{bmatrix} k-j-2 \\ j \end{bmatrix}_{q^2} - q^{2k-2-4j}\begin{bmatrix} k-j-2 \\ j-1 \end{bmatrix}_{q^2}$$
$$+ q^{k-1-2j}\left( \begin{bmatrix} k-j-1 \\ j \end{bmatrix}_{q^2} - q^{2j}\begin{bmatrix} k-j-2 \\ j \end{bmatrix}_{q^2} - \begin{bmatrix} k-j-2 \\ j-1 \end{bmatrix}_{q^2} \right) = 0$$

by using both recurrences for the $q$-binomial coefficients.

From
$$\frac{r_{2n+1}(1,-q^{k+1})}{r_{2n}(1,-1)} - \frac{r_{2n+1}(1,-q^k)}{r_{2n}(1,-1)} = q^k(1-q^{2n+1})\frac{r_{2n}(1,-q^k)}{r_{2n}(1,-1)}.$$

we see that
$$c(n,k) = \frac{r_{2n}(1,-q^k)}{r_{2n}(1,-1)} = q^{-k}(b(n+1,k+1) - b(n+1,k)). \tag{2.13}$$

This gives

$$c(n,k) = \sum_{j=0}^{\lfloor k/2 \rfloor} (-1)^j q^{j^2} q^{2jn} \frac{[k]}{[2k-2j]} \begin{bmatrix} k-j \\ j \end{bmatrix}_{q^2} \prod_{i=1}^{k-2j} (1+q^i).$$

Therefore we have proved

**Theorem 2**

*For each fixed $k \geq 0$ the following identities hold:*

$$\frac{r_{2n-1}(1,-q^k)}{r_{2n}(1,-1)} = \sum_{j=0}^{\lfloor (k-1)/2 \rfloor} (-1)^j q^{j^2} q^{2jn} \begin{bmatrix} k-j-1 \\ j \end{bmatrix}_{q^2} \prod_{i=1}^{k-1-2j} (1+q^i) \qquad (2.14)$$

*and*

$$\frac{r_{2n}(1,-q^k)}{r_{2n}(1,-1)} = \sum_{j=0}^{\lfloor k/2 \rfloor} (-1)^j q^{j^2} q^{2jn} \frac{[k]}{[2k-2j]} \begin{bmatrix} k-j \\ j \end{bmatrix}_{q^2} \prod_{i=1}^{k-2j} (1+q^i). \qquad (2.15)$$